\input amsppt.sty
\documentstyle{amsppt}
\font\tenrm=cmr10
\hsize=4.75in
\vsize=8in
\NoBlackBoxes

\font\bb=msbm10 at 10 pt
\font\bbb=msbm10 at 8 pt
\define\op#1{\operatorname{\fam=0\tenrm{#1}}} 
\def\croi{\hbox{$\mathop{\mathrel\times\joinrel\mathrel{\vrule height 5pt
 depth 0pt}}$}}

\def\C{\hbox{\bb C}}

\def\N{\hbox{\bb N}}
\def\F{\hbox{\bb F}}
\def\n{\hbox{\bbb N}}

\def \e {\hbox{End}}
\def \p {\hbox{Prob}}


\rightheadtext {Amenability and exactness  }
\leftheadtext {Claire Anantharaman-Delaroche}
\topmatter

\title
Amenability and exactness for dynamical systems and their
$C^*$-algebras
\endtitle

\author
Claire Anantharaman-Delaroche
\endauthor

\affil
UMR 6628, MAPMO\\
Universit\'e d'Orl\'eans-CNRS
\endaffil

\address
D\'epartement de Math\'ematiques, Universit\'e d'Orl\'eans,
B. P. 6759, F-45067 Orl\'eans Cedex 2
\endaddress

\email
claire\@labomath.univ-orleans.fr
\endemail

\date
April 2000
\enddate

\abstract
In this survey, we study the relations between amenability (resp.
amenability at infinity)
of $C^*$-dynamical systems and equality or nuclearity (resp. exactness) of the
corresponding crossed products.
\endabstract
\endtopmatter

\document

\heading
1. Introduction
\endheading

Recently, amenable transformation groups have attracted much attention in
connection
with the Novikov conjecture for discrete groups, following the works of Yu,
Higson-Roe and
Higson (\cite{\bf Yu},\cite{\bf HR},\cite{\bf Hi}). Therefore I was
encouraged to write down
my simplified version of my paper \cite{\bf AD},
where I studied the relations between amenability of $C^*$-dynamical
systems and nuclearity, as well as equality, of the associated full and
reduced
crossed products.

These results are also contained in the monograph \cite{\bf ADR}
that Jean Renault and I have
recently written on amenable groupoids. However, the proofs given in
Sections 2 and 3 below
are intended to be more direct and more easily accessible to readers
unfamiliar with
groupoids. Another advantage is that, when working with transformation
groups, one does not
need any separability assumption on the locally compact spaces.

In addition, as an immediate consequence of the techniques developed here,
it is shown in
Section 4 that a discrete group is exact if and only if it admits an
amenable action
on a compact space. For a nice self-contained proof of this result we refer to
 \cite{\bf Oz}.

These notes are part of a series of talks that I gave for the Instructional
Conference in Operator Algebras and Operator Spaces in Edinburgh during
the first two weeks of April 2000. The main results are Theorems 3.4 and 4.6.

\heading
2. Amenable transformation groups
\endheading

A {\it transformation group} is a left $G$-space $X$, where $X$ is a
locally compact space,
$G$ is a locally compact group and $(x,s) \mapsto s.x$ is a continuous left
action
 from $X \times G$
to $X$. Let us first introduce  some notations. We denote by $C_c(G)$ the space
of complex valued continuous functions with compact support on $G$, and
$\p(G)$ will be
the set of probability measures on $G$, equipped with the weak*-topology.
Given $f\in C_c(G)$,
$x\in G$, $m\in \p(G)$ and $s\in G$, we set
$$(s.f)(x) = f(s^{-1}.x), \quad (s.m)(f) = m(s^{-1}.f).$$

If $g\in C_c(X\times G)$, then $g^x$ will be the map $t\mapsto g^x(t) =
g(x,t)$,
and $g(t)$ will be the map $x\mapsto g(t)(x) = g(x,t)$. Finally, $dt$
denotes the
left Haar measure on $G$.

\definition {2.1 Definition} {\rm We say that the transformation group
$(X,G)$ (or
that the $G$-action on $X$) is {\it amenable} if there exists a net
$(m_i)_{i\in I}$
of continuous maps $x \mapsto m^{x}_i$ from $X$ into the space $\p(G)$ such
that
$$\lim_i \| s.m^{x}_i - m_{i}^{s.x}\|_1 = 0$$ uniformly on compact subsets of
$X\times G$.

Such a net $(m_i)_{i\in I}$ will be called an {\it approximate invariant
continuous mean
} (a.i.c.m. for short).
\enddefinition

Before giving examples, let us state several equivalent definitions.

\proclaim{2.2 Proposition} The following conditions are equivalent:
\roster
\item $(X,G)$ is an amenable transformation group.
\item There exists a net $(g_i)_{i\in I}$ of nonnegative continuous functions
on
$X\times G$ such that
\itemitem {\rm a)} for every $i\in I$ and $x\in X$, $\int_G g_{i}^x(t)dt = 1$;
\itemitem{\rm b)} $\lim_i \int_G |g_{i}^{s.x}(st) - g_{i}^{x}(t)|dt =0$
uniformly on compact
subsets of $X\times G$.
\item There exists a net $(g_i)_{i\in I}$ in $C_c(X\times G)^+$ such that
\itemitem {\rm a)} $\lim_i \int_G g_{i}^{x} dt = 1$ uniformly on compact
subsets of $X$;
\itemitem {\rm b)} $\lim_i \int_G |g_{i}^{s.x}(st) - g_{i}^{x}(t)|dt =0$
uniformly on compact
subsets of $X\times G$.
\endroster
\endproclaim

\demo{Proof}
(1) $\Rightarrow$ (2). Let $g\in C_c(G)^+$ such that $\int g(t) dt =1$ and set
$$g_i(x,s) = \int g(t^{-1}s) dm_{i}^{x}(t).$$
By the Fubini theorem, we get $\int g_{i}^x(t) dt = 1$ for every $x\in X$
and $i\in I$.
Moreover, for $(x,s)\in X\times G$ we have
$$
\aligned
\int |g_{i}^{s.x}(st) - g_{i}^{x}(t)| dt
&= \int|\int g(u^{-1}st) dm_{i}^{s.x}(u)
-\int g(u^{-1}t) dm_{i}^{x}(u)| dt\\
&= \int | \int g(u^{-1}st) dm_{i}^{s.x}(u) -\int g(u^{-1}st) d
(s.m_{i}^{x})(u)|dt\\
&\leq \int\int g(u^{-1}st) d|m_{i}^{s.x} - s.m_{i}^{x}|(u) dt
\endaligned
$$
where $|m_{i}^{s.x} - s.m_{i}^{x}|$ is the total variation of $m_{i}^{s.x}
- s.m_{i}^{x}$.
Using again the Fubini theorem we obtain the majoration
$$\leq \int [\int g(u^{-1}st) dt] d|m_{i}^{s.x} - s.m_{i}^{x}|(u) \leq
\| m_{i}^{s.x} - s.m_{i}^{x}\|_1.$$
This tends to $0$ uniformly on compact subsets of $X\times G$.

(2) $\Rightarrow$ (3). An easy approximation argument allows to replace the net
$(g_i)_{i\in I}$ by a net in $C_c(X\times G)^+$ satisfying conditions
a) and b) of (3).

(3) $\Rightarrow$ (2). Let $(g_i)_{i\in I}$ as in (3) and choose $g\in
C_c(G)^+$ as above.
Let us define $(g_{i,n})$ by
$$g_{i,n}(x,s) = \frac{g_{i}(x,s) + \frac{1}{n}g(s)}{\int g_i(x,t)dt +
\frac{1}{n}}.$$
Then $(g_{i,n})_{(i,n) \in I\times \n^*}$ satisfies conditions a) and b) of
(2).

(2) $\Rightarrow$ (1) is obvious: given $(g_i)_{i\in I}$ as in (2), we define
$m_{i}^x$ to be the probability measure with density $g_{i}^x$ with respect
to the left Haar
measure. Then $(m_i)_{i\in I}$ is an a.i.c.m.
\qed
\enddemo

\definition {2.3 Definition} {\rm A function $h$ defined on $X\times G$ is
said to be of
{\it positive type} if for every $x\in X$, every positive integer $n$ and
every $t_1,\dots, t_n \in G$, $\lambda_1,\dots,\lambda_n \in \C$ we have
$$\sum_{i,j} h(t_{i}^{-1}.x, t_{i}^{-1}t_j) \overline{\lambda}_i \lambda_j
\geq 0.$$}
\enddefinition

We denote by $L^2(G)\otimes C_0(X)$ the completion of $C_c(X\times G)$
endowed with the norm
$$\|\xi\|_2 = (\sup_{x\in X} \int|\xi(x,t)|^2 dt)^{1/2}.$$
Remark that $L^2(G)\otimes C_0(X)$ is a Hilbert $C^*$-module over the algebra
$C_0(X)$ of continuous functions vanishing at infinity, with scalar product
$$\langle \xi,\eta\rangle(x) = \int \overline{\xi(x,t)} \eta(x,t) dt,$$
and module structure defined by
$$(\xi f)(x,s) = \xi(x,s) f(x).$$
For basic facts on Hilbert $C^*$-modules we refer to \cite{\bf La}.

To every $\xi \in L^2(G)\otimes C_0(X)$, we associate the coefficient
$(\xi,\xi)$
defined by
$$(\xi,\xi)(x,s) =  \int \overline{\xi(x,t)} \xi(s^{-1}.x, s^{-1}t) dt.$$
It is a continuous positive type function on $X\times G$, with compact
support if $\xi$ is so. The
following result is the analogue for transformation groups of a well known
result of
Godement \cite{{\bf Di}, Th\'eor\`eme 13.8.6}.

\proclaim{2.4 Lemma} Let $k$ be a continuous function of
positive type on $X\times G$, with compact support. There exists $\xi \in
L^2(G)\otimes C_0(X)$ such that $k = (\xi,\xi)$.
\endproclaim

\demo{Proof} Let $\rho(k)$ be the $C_0(X)$-linear bounded endomorphism of
$L^2(G)\otimes C_0(X)$ such that
$$(\rho(k)\xi)(x,s) = \int k(s^{-1}.x, s^{-1}t)\xi(x,t) dt$$
for every $\xi\in L^2(G)\otimes C_0(X)$. Since $k$ is a positive type
function,
$\rho(k)$ is a positive operator. Let $(f_i)_{i\in I}$ be an approximate
unit of
$C_0(X)$ and $(\varphi_j)_{j\in J}$ be a standard approximate unit of
$L^1(G)$ (i.e.
made of nonnegative continuous functions with compact supports in a basis
of neighbourhoods
of the unit of $G$). We set
$$\xi_{i,j} = \rho(k)^{1/2}(f_i\otimes\varphi_j).$$
Then $(\xi_{i,j})$ is a Cauchy net in $L^2(G)\otimes C_0(X)$, and if $\xi$
denotes
its limit, we have
$$\aligned
(\xi,\xi)(x,s) &= \lim_{i,j} (\xi_{i,j},\xi_{i,j})(x,s)\\
&= \lim_{i,j} f_i(x) f_i(s^{-1}.x) \int\int \varphi_j(t)\varphi_j(r)
k(t^{-1}.x, t^{-1}sr)dt dr\\
&= k(x,s).
\endaligned$$
\qed
\enddemo

\proclaim{2.5 Proposition} The following conditions are equivalent:
\roster
\item $(X,G)$ is an amenable transformation group.
\item There exists a net $(\xi_i)_{i\in I}$ in $C_c(X\times G)$ such that
\itemitem{\rm a)} $\lim_i \int |\xi_i(x,t)|^2 dt = 1$ uniformly on compact
subsets
of $X$;
\itemitem{\rm b)} $\lim_i \int|\xi_i(s.x,st) - \xi_i(x,t)|^2 dt = 0$
uniformly on compact subsets of $X\times G$.
\item There exists a net $(h_i)_{i\in I}$ of positive type functions in
$C_c(X\times G)$ which tends to $1$ uniformly on compact subsets of
$X\times G$.
\endroster
\endproclaim

\demo{Proof} (1) $\Rightarrow$ (3). Let $(g_i)_{i\in I}$ be a net in
$C_c(X\times G)^+$ which
fulfills conditions (a) and (b) of (3) in Proposition 2.2. Wet set $\xi_i =
\sqrt{g_i}$
and $h_i = (\xi_i,\xi_i)$. Then
$$h_i(s^{-1}.x,e) + h_i(x,e) - 2h_i(x,s) =
\int|\xi_i(s^{-1}.x,s^{-1}t) - \xi_i(x,t)|^2 dt. \quad\quad (*)$$
The inequality $|a - b|^2 \leq |a^2 - b^2|$ for $a,b \geq 0$ gives the
majoration
$$\int |\xi_i(s^{-1}.x,s^{-1}t) - \xi_i(x,t)|^2 dt \leq \int
|g_i(s^{-1}.x,s^{-1}t) - g_i(x,t)| dt.
\quad \quad (**)$$
Since $h_i(x,e) = \int g_i(x,t)dt$, we see from $(*)$ and $(**)$ that
$(h_i)_{i\in I}$
goes to $1$ uniformly on compact subsets of $X\times G$.

The implication (3) $\Rightarrow$ (2) follows immediately from Lemma 2.4
and $(*)$. We just
have to observe that the elements of $L^2(G)\otimes C_0(X)$ can be
approximated by
elements of $C_c(X\times G)$.

(2) $\Rightarrow$ (1). Let $(\xi_i)_{i\in I}$ as in (2) and set $g_i =
|\xi_i |^2$
for $i\in I$. Using the Cauchy-Schwarz inequality and the inequality
$$| |a|^2 - |b|^2| \leq (|a| + |b|) |a - b|,$$
we get
$$\aligned
\int &|g_i(s.x,st) - g_i(x,u)|dt \\
&\leq [\int (|\xi_i(s.x,st)| + |\xi_i(x,t)|)^2 dt]^{1/2}
\times [\int |\xi_i(s.x,st) - \xi_i(x,t)|^2 dt]^{1/2}\\
&\leq [(\int |\xi_i(s.x,t)|^2)^{1/2} + (\int |\xi_i(x,t)|^2 dt)^{1/2}]
\times [\int |\xi_i(s.x,st) - \xi_i(x,t)|^2 dt]^{1/2}.
\endaligned
$$
This last term goes to $0$ uniformly on compact subsets of $X\times G$.
\qed
\enddemo

\remark{2.6 Remark}  Obviously, when $X$ and $G$ are $\sigma$-compact, we may
replace nets by sequences in the above statements.
\endremark

\example{2.7 Examples} (1) If $X$ is reduced to a point, our definition of
amenability,
under the form given Proposition 2.2, is the existence of a net
$(k_i)_{i\in I}$
of continuous nonnegative functions such that $\lim_i \int k_i(t) dt = 1$ and
$s\mapsto \int |k_i(st) - k_i(t)| dt$ goes to $0$ uniformly on compact subsets
of $G$; it is essentially Reiter's condition $(P_1)$ \cite{{\bf Gre}, page
44}. The equivalent
condition given in Proposition 2.5 (3), namely the existence of a net
$(h_i)_{i\in I}$
 of continuous positive type compactly supported functions converging to $1$
uniformly on compact subsets of $G$ is essentially the weak containment
condition
introduced by Godement \cite{{\bf Gre}, page 61}.

(2) If $G$ is an amenable group, every transformation group $(X,G)$
is amenable. Indeed, given $(k_i)_{i\in I}$
as in (1) above, we define $m_i$ to be the constant map on $X$ whose value
is the probability
measure with density $k_i$ with respect to the left Haar measure. Then
$(m_i)_{i\in I}$ is an
a.i.c.m.

Note that if $(X,G)$ is amenable and if $X$ has a $G$-invariant probability
measure $\mu$,
then $G$ is an amenable group. Indeed, let us consider a net $(g_i)_{i\in
I}$ as in Proposition
2.2 (3), and set $k_i(t) = \int  g_i(x,t) d\mu(x)$. Then $(k_i)_{i\in I}$
has the properties stated
in (1), and therefore $G$ is amenable.

(3) For every locally compact group $G$, its left action on itself is
amenable. Indeed, let us
choose a nonnegative function $g \in C_c(G)^+$, such that $\int g(t) dt =
1$ and set
$m^x(f) = \int f(t) g(x^{-1}t) dt$ for every $x\in G$ and every $f\in C_c(G)$.
Then $x\mapsto m^x$ is continuous from $G$ into $\p(G)$ and
$m^{sx} = s.m^x$ for every $(x,s) \in G\times G$.

Such an action, having a continuous
invariant system of probability measures, is called {\it proper}. This is
equivalent to
the usual properness of the map $(x,t)\mapsto (x, t.x)$ from $X\times G$
into $X\times X$
(see \cite{{\bf ADR}, Cor. 2.1.17}).

(4) Let $\F_2$ be the free group with two generators $a$ and $b$. The boundary
$\partial\F_2$ is the set of all infinite reduced words $\omega = a_1 a_2
\dots a_n \dots$
in the alphabet $S = \{a, a^{-1}, b, b^{-1}\}$. There is a natural topology
on $\partial\F_2$ for which $\partial\F_2$ is the Cantor discontinuum, and
$\F_2$ acts
continuously to the left by concatenation. The transformation group
$(\partial\F_2,\F_2)$
is amenable. Indeed, for $n\geq 1$ and $\omega = a_1 a_2 \dots a_n \dots$,
define
$$\aligned
m_{n}^\omega(\{t\}) &= \frac{1}{n} \quad\hbox{if} \quad t = a_1\dots a_k,
\quad k\leq \frac{1}{n},\\
&= 0 \quad\hbox{otherwise}.
\endaligned$$
Then $(m_n)_{n\geq 1}$ is an a.i.c.m.

More generally, every hyperbolic group $\Gamma$ acts amenably on its Gromov
boundary
$\partial\Gamma$.
This was first proved by Adams \cite{\bf Ad}, by showing that for every
quasi-invariant measure
$\mu$ on $\partial\Gamma$, the measured transformation group
$(\partial\Gamma,\Gamma,\mu)$
is amenable in the sense of Zimmer. This is equivalent to the (topological)
amenability
of $(\partial\Gamma,\Gamma)$ defined in these notes (see \cite{{\bf ADR},
Th. 3.3.7}). More
recently, Germain has given a direct proof by constructing explicitly an
approximate
invariant mean \cite{{\bf ADR}, Appendix B}.

(5) Another convenient way to show that transformation groups are amenable
is to use the
invariance of this notion by Morita equivalence \cite{{\bf ADR}, Th.
2.2.17}. Let us consider for
instance a locally compact group $G$ and two closed subgroups $H$ and $K$.
Then the
left $H$-action on $G/K$ and the right $K$-action on $H\setminus G$ are Morita
equivalent. If $K$ is an amenable group, it follows that $(G/K,H)$ is an
amenable transformation group.
\endexample

\heading
3. Amenability and crossed products
\endheading

Let us recall first a few facts on crossed products. We refer to \cite{{\bf
Pe}, Chap. 7} for more details.

\definition{3.1 Definitions}  Let $G$ be a locally compact group. A {\it
$G$-$C^*$-algebra}
is a $C^*$-algebra endowed with a continuous action $\alpha$ of $G$.
Specifically,
$\alpha$ is a homomorphism from $G$ to the group of automorphisms of $A$,
such that for $a \in A$
the map $s\mapsto \alpha_s(a)$ is norm continuous. We shall often write $s.a$
for $\alpha_s(a)$.

A {\it covariant representation} of the $G$-$C^*$-algebra $A$ is a pair
$(\pi,\sigma)$
where $\pi$ and $\sigma$ are representations of $A$ and $G$ respectively in
a Hilbert space
$H$, such that
$$\sigma(s)\pi(a)\sigma(s)^{-1} = \pi(s.a)$$
for every $a\in A, s\in G$.
\enddefinition

We denote by $C_c(G,A)$ the space of continuous functions with compact
support from $G$ to
$A$, equipped with the product
$$f*g(s) = \int f(t)\alpha_t(g(t^{-1}s)) dt$$
and the involution
$$f^*(s) = \frac{\alpha_s(f(s^{-1}))^*}{\Delta(s)}$$
where $\Delta$ is the modular function of $G$.

Given a covariant representation $(\pi,\sigma)$, we set, for $f\in C_c(G,A)$,
$$(\pi\times\sigma)(f) = \int \pi(f(t))\sigma(t) dt \in {\Cal B}(H).$$
Then $\pi\times\sigma$ is a representation of the involutive algebra
$C_c(G,A)$. The {\it full crossed product} $C^*(G,A)$ is the completion of
$C_c(G,A)$ with respect to the norm
$$\| f \| = \sup_{(\pi,\sigma)} \| (\pi\times \sigma) (f)\|$$
where $(\pi,\sigma)$ ranges over all covariant representations of $A$.

Note that $\pi\times \sigma$ extends to a representation of $C^*(G,A)$
and that every non-degenerate representation of $C^*(G,A)$ arises in this way.

We now view the vector space $C_c(G,A)$ as a right $A$-module
in the obvious way and define on it the $A$-valued scalar product
$$\langle f,g\rangle_A = \int f(t)^*g(t) dt.$$
As we did previously for $A = C_0(X)$, we denote by $L^2(G)\otimes A$ the
Hilbert $C^*$-module over $A$ obtained
from $C_c(G,A)$ by completion. The $C^*$-algebra of $A$-linear
bounded operators on $L^2(G)\otimes A$ admitting an adjoint will be denoted
by $\e_A (L^2(G)\otimes A)$. We define an involutive homomorphism
$$\Lambda : C_c(G,A) \rightarrow \e_A (L^2(G)\otimes A)$$
by the formula
$$\Lambda(f) \xi(s) = \int\alpha_{s^{-1}}(f(t))\xi(t^{-1}s) dt,$$
and the {\it reduced crossed product} $C^*_{r}(G,A)$ is the closure of
$\Lambda(C_c(G,A))$.

Given a non-degenerate representation $\pi : A \rightarrow {\Cal B}(H)$,
the Hilbert space $L^2(G,A) \otimes_A H$ will be identified to $L^2(G,H)$
by the map sending $f\otimes \xi$ to $t \mapsto \pi(f(t))\xi$. Denote by
$\lambda$
the left regular representation of $G$ into $L^2(G)$ and by $\tilde{\pi}$ the
representation of $A$ into $L^2(G,H)$ such that
$$\tilde{\pi}(a)\xi(t) = \pi(\alpha_{t^{-1}}(a))\xi(t), \quad \forall a\in
A, \forall t\in G.$$
Then $(\tilde{\pi}, \lambda\otimes 1_H)$ is a covariant representation of
the $G$-$C^*$-algebra
$A$. We denote by $\tilde{\pi}\times \lambda$ (instead of
$\tilde{\pi}\times (\lambda\otimes 1_H)$)
the corresponding representation of $C^*(G,A)$, said to be {\it induced by}
$\pi$. If
$\pi$ is a faithful representation, then $b \mapsto \Lambda(b) \otimes_A 1_H$
is a faithful representation of $C^*_{r}(G,A)$ which coincides with
$\tilde{\pi}\times \lambda$
on $C_c(G,A)$. It follows that $C^*_{r}(G,A)$ is a quotient of $C^*(G,A)$
and that
every induced representation passes to the quotient.

Let $(X,G)$ be a transformation group; the $C^*$-algebra $C_0(X)$ of
continuous functions on $X$ vanishing at
infinity is a $G$-$C^*$-algebra, with
$s.f(x) = f(s^{-1}.x)$ for $s\in G, x\in X, f\in C_0(X)$. The corresponding
crossed products
will be denoted by $C^*(X\croi \,G)$ and $C^*_{r}(X\croi \,G)$ respectively.
When $X$ is reduced to a point, we get the
{\it full} and {\it reduced} $C^*$-algebras of $G$, denoted by $C^*(G)$ and
$C_{r}^*(G)$ respectively.

If $G$ is a discrete group, note that $C_0(X) \subset C_{r}^*(X\croi\, G)$.
We denote by
$\lambda(s)$ every $s\in G$ viewed as a unitary of $C_{r}^*(G)$ and also as
an element
of the multiplier algebra $M(C_{r}^*(X\croi\, G))$. Then for $\varphi \in
C_0(X)$,
$\lambda(s) * \varphi$ and $\varphi *\lambda(s)$ will denote the obvious
products in
$M(C_{r}^*(X\croi\, G))$, and belong to $C_{r}^*(X\croi\, G)$.

In the proof of our main theorem 3.4, we  shall need the following lemma.

\proclaim{3.2 Lemma} Let $G$ be a discrete group, and let $P :
C_{r}^*(X\croi\,G) \rightarrow
M_n(\C)$ be a completely positive map. Given $\varepsilon > 0$, $\varphi
\in C_0(X)$
and a finite subset $F \subset G$, there exists a completely positive map
$\tilde{P} : C_{r}^*(X\croi\,G) \rightarrow M_n(\C)$ such that
\roster
\item $\| P(\varphi * \lambda(s)*\varphi) - \tilde{P}(\varphi *
\lambda(s)*\varphi)
\| \leq \varepsilon$ for every $s\in F$;
\item $s \mapsto \tilde{P}(\phi * \lambda(s)*\psi)$ has a finite support
for every
$\phi,\psi \in C_0(X)$.
\endroster
\endproclaim

\demo{Proof} Since $P$ is a completely positive map, by the Stinespring
dilatation theorem
there exist a representation $\rho$ of $C_{r}^*(X\croi\, G)$ into
$H_{\rho}$ and
vectors $e_1,\dots,e_n \in H_{\rho}$ such that for $a \in C_{r}^*(X\croi\, G)$,
$$P(a) = [\langle e_i, \rho(a) e_j\rangle] \in M_n(\C).$$
Let us choose a faithful representation $\pi$ of $C_0(X)$ in a Hilbert
space $H$.
Since $\tilde{\pi}\times \lambda$ is a faithful representation of
$C_{r}^*(X\croi \,G)$
into $l^2(G,H)$, given $\eta > 0$ there exist a multiple $(\tilde{\pi}\times \lambda)
\otimes 1_K$ of this representation (that is $\tilde{\pi}_K \
times \lambda_K$,
where $\pi_K = \pi\otimes 1_K$ and $\lambda_K = \lambda \otimes 1_{H\otimes
K}$),
and vectors $\xi_1,\dots,\xi_n$ in $l^2(G, H\otimes K)$ such that
$$|\langle e_i, \rho(\varphi * \lambda(s)*\varphi) e_j\rangle - \langle \xi_i,
(\tilde{\pi}_K \times \lambda_K)(\varphi * \lambda(s)*\varphi
)\xi_j\rangle| < \eta$$
for $s\in F$, $i,j \in \{1,\dots,n\}$. Moreover the $\xi_i$'s may be chosen
with
finite support. Since
$$\langle \xi_i, (\tilde{\pi}_K \times \lambda_K)(\phi * \lambda(s)*\psi
)\xi_j\rangle = \langle \tilde{\pi}_K(\overline{\phi})\xi_i,
\lambda_K(s)\tilde{\pi}_K(\psi)\xi_j
\rangle$$
where $\tilde{\pi}_K(\overline{\phi})\xi_i$ and $\tilde{\pi}_K(\psi)\xi_j$
have finite support,
it follows that the map $$a \mapsto \tilde{P}(a) = [\langle \xi_i,
(\tilde{\pi}_K \times \lambda_K)(a)
\xi_j \rangle]$$
satisfies the required properties, provided $\eta$ is small enough. \qed
\enddemo

\definition{3.3 Definitions} Let $X$ be a locally compact space. A {\it
$C_0(X)$-algebra}
is a $C^*$-algebra $A$ equipped with a morphism from $C_0(X)$ into the
center $Z(M(A))$ of the multiplier
algebra of $A$, such that $\overline{C_0(X)A} = A.$

Let $(X,G)$ be a transformation group. A {\it $G$-$C_0(X)$-algebra} is a
$G$-$C^*$-algebra
$A$ which is also a $C_0(X)$-algebra and satisfies the compatibility condition
$$s.(fa) = (s.f) (s.a), \quad \forall f\in C_0(X), a\in A, s \in G.$$
\enddefinition

\proclaim{3.4 Theorem \rm\cite{\bf AD}} Let $(X,G)$ be a transformation
group, and consider the
following conditions :
\roster
\item $(X,G)$ is amenable.
\item For every $G$-$C_0(X)$-algebra $A$, $C_{r}^*(G,A) = C^*(G,A)$.
\item For every nuclear $G$-$C_0(X)$-algebra $A$, $C_{r}^*(G,A)$ is nuclear.
\item $C_{r}^*(X\croi\,G)$ is nuclear.
\endroster
Then (1) $\Rightarrow$ (2) $\Rightarrow$ (3) $\Rightarrow$ (4). Morever
(4) $\Rightarrow$ (1) if $G$ is discrete.
\endproclaim

\demo{Proof} (1) $\Rightarrow$ (2). Let $(\pi,\sigma)$ be a covariant
representation
of the $G$-$C^*$-algebra $A$ in $H$. We show that $\pi\times \sigma$ is
weakly contained in the
representation $\tilde{\pi}\times \lambda$. Let $(\xi_i)$ be a net as in
Proposition 2.5 (2). Given $\xi \in H$, we define, for $t\in G$,
$$\eta_i(t) = \sigma(t^{-1})\pi(\xi_i(t))\xi.$$
Here $\xi_i(t)\in C_c(X)$ is viewed as an element of $Z(M(A))$, and $\pi$
is extended to $M(A)$. For $f\in C_c(X\times G)$, a straightforward computation
gives
$$\aligned
\langle \eta_i, (\tilde{\pi}\times\lambda)(f) \eta_i\rangle &=
\int\langle \eta_i(s), \pi(\alpha_{s^{-1}}(f(t))\eta_i(t^{-1}s)\rangle ds dt\\
&= \int \langle \pi(h_i(t))\xi,\pi(f(t))\sigma(t)\xi\rangle dt
\endaligned$$
where $h_i(x,s) = \int \overline{\xi_i(x,t)}\xi_i(s^{-1}.x, s^{-1}t) dt$.
Since $(h_i)$ goes to $1$ uniformly on compact subsets of $X\times G$,
we have
$$\lim_i \langle \eta_i, (\tilde{\pi}\times\lambda)(f) \eta_i\rangle = \langle
\xi, (\pi\times \sigma)(f) \xi\rangle.$$
Therefore $\pi\times \sigma$ is weakly contained in $\tilde{\pi}\times
\lambda$,
and (2) is proved.

(2) $\Rightarrow$ (3). Let $B$ be a $C^*$-algebra. Since $A\otimes B$ is a
$G$-$C_0(X)$-algebra in an obvious way, we have
$$C_{r}^*(G,A)\otimes B = C_{r}^*(G,A\otimes B) = C^*(G,A\otimes B)$$
by (2). Moreover $A\otimes B = A\otimes_{max}B$ because $A$ is nuclear. It
follows
that
$$\aligned
C_{r}^*(G,A)\otimes B = C^*(G,A\otimes_{max} B) &=C^*(G,A)\otimes_{max} B\\
& = C_{r}^*(G,A)\otimes_{max} B,
\endaligned$$
and therefore $C_{r}^*(G,A)$ is nuclear.

(3) $\Rightarrow$ (4) is obvious.

Assume now that $G$ is discrete and let us show that (4) $\Rightarrow$ (1).
Given $\varepsilon > 0$ and a compact subset $K\times F$ of $X\times G$, we
will
prove the existence of a compactly supported  positive type function $h$ such
that $$\displaystyle \sup_{(x,s)\in K\times F} |h(x,s) - 1| \leq \varepsilon.$$
Let $\varphi$ be a continuous function with compact support on $X$, with values
in $[0,1]$ such that $\varphi(x) = 1$ for $x\in K\cup F^{-1}K$.
Given $s \in G$, define $\varphi_s$ by
$$\aligned
\varphi_s(x,t) &= 0 \quad\hbox{if}\quad t\not= s, \\
 \varphi_s(x,s) &= \varphi(s.x).
\endaligned$$
We view $\varphi_s$ as an element of the Hilbert $C^*$-module
$l^2(G)\otimes C_0(X)$.

Now consider a completely positive contraction $\Phi : C_{r}^*(X\croi\, G)
\rightarrow C_{r}^*(X\croi\, G)$ which factorizes through a matrix algebra :
$\Phi = Q\circ P$ where $P : C_{r}^*(X\croi\, G) \rightarrow M_n(\C)$ and
$Q : M_n(\C) \rightarrow C_{r}^*(X\croi\, G)$ are
completely positive contractions. Define
$$h(x,s) = \langle \varphi_s, \Phi(\varphi *\lambda(s) *\varphi) \varphi_e
\rangle
(s^{-1}.x).$$
A straightforward computation shows that
$$h(s^{-1}.x, s^{-1}t) = \langle \varphi_{s^{-1}}, \Phi(\varphi
*\lambda(s^{-1}t)
*\varphi) \varphi_{t^{-1}} \rangle(x),$$
and therefore $h$ is a positive type function since $\Phi$ is completely
positive.

Since $C_{r}^*(X\croi\, G)$ is nuclear, we may choose $\Phi$ such that
$$\| \Phi(\varphi *\lambda(s) *\varphi) - \varphi *\lambda(s) *\varphi\|
\leq \varepsilon/2$$
for every $s \in F$. Then
$$\sup_{x\in X} |h(x,s) - \langle \varphi_s, \Lambda(\varphi *\lambda(s)
*\varphi)
\varphi_e \rangle(s^{-1}.x)| \leq \varepsilon/2$$ for $s\in F$.
Now we have
$$\langle \varphi_s, \Lambda(\varphi *\lambda(s) *\varphi)
\varphi_e \rangle(s^{-1}.x) = \varphi(x)^2\varphi(s^{-1}.x)^2 = 1$$
for $(x,s) \in K\times F$. The only point is that $h$ needs not have a
compact support, but using Lemma
3.2, we construct a compactly supported positive type function $\tilde{h}$
such that $\displaystyle \sup_{(x,s)\in X\times F} |\tilde{h}(x,s) -
h(x,s)| \leq
\varepsilon/2$. Indeed, we choose a completely positive map $\tilde{P} :
 C_{r}^*(X\croi\, G) \rightarrow M_n(\C)$ such that
$$\| P(\varphi * \lambda(s)*\varphi) - \tilde{P}(\varphi * \lambda(s)*\varphi)
\| \leq \varepsilon/2$$ for every $s\in F$ and such that
$s \mapsto \tilde{P}(\varphi * \lambda(s)*\varphi)$ has a finite support. Then
we consider $\tilde{\Phi} = Q \circ \tilde{P}$ and we set
$$\tilde{h}(x,s) = \langle \varphi_s, \tilde{\Phi}(\varphi *\lambda(s)
*\varphi) \varphi_e \rangle
(s^{-1}.x).$$
Note that $\tilde{h}(x,s) = \varphi(x)\tilde{\Phi}(\varphi *\lambda(s)
*\varphi) \varphi_e \rangle
(s^{-1}.x,s)$ has already its support in $(\hbox{Supp}\, \varphi) \times
G$. \qed

\enddemo

\remark{3.5 Remark} The implication (4) $\Rightarrow$ (1) is false when $G$
is not discrete.
Indeed, Connes proved in \cite{\bf Co} that $C^*(G)$ is nuclear for every
connected
locally compact separable group.

On the other hand, for every locally compact group $G$, we have
$C_{r}^*(G) = C^*(G)$ if and only if $G$ is amenable \cite{{\bf Pe}, Th.
7.3.9}. We do not know, even for a
discrete group $G$, whether $C_{r}^*(X\croi\, G) = C^*(X\croi\, G)$
implies the amenability of the transformation group.
\endremark

\heading
4. Exactness and amenability at infinity
\endheading

\definition{4.1 Definition} We say that a locally compact group $G$ is
{\it amenable at infinity} if it admits an amenable action on a compact
space $X$.
\enddefinition

\example{4.2 Examples} (1) An amenable locally compact group is amenable at
infinity:
one can take $X$ reduced to a point.

(2) A discrete group $G$ which is hyperbolic in the sense of Gromov is
amenable at infinity
since its action on its Gromov boundary is amenable.

(3) Every closed subgroup $H$ of a connected Lie group $G$ is amenable at
infinity.
Indeed, $G$ possesses a closed amenable subgroup $K$ such that $G/K$ is
compact. Then
$(G/K, H)$ is an amenable transformation group.
\endexample

Let us assume now that $G$ is a countable discrete group. We give a
characterization
of amenability at infinity of $G$ which only invokes $G$. A set $T_F$ of
the form
$\{(s,t)\in G\times G : s^{-1}t \in F\}$, where $F$ is a finite subset of $G$,
will be called a {\it tube}.

\proclaim{4.3 Proposition}  Let $G$ be a countable discrete group. The
following
conditions are equivalent:
\roster
\item $G$ is amenable at infinity.
\item The action of $G$ on its Stone-\u{C}ech compactification, extending
the left multiplication,
is amenable.
\item $C_{r}^*(\beta G\croi\, G)$ is nuclear.
\endroster
\endproclaim

\demo{Proof} Of course we have (3) $\Leftrightarrow$ (2) $\Rightarrow$ (1).
Now,
assume that there exists an amenable transformation group $(X,G)$ with $X$
compact. Choosing $x_0 \in X$, the map $s\mapsto s.x_0$ from $G$ to $X$ extends
to a continuous map $p: \beta G \rightarrow X$, by the universal property
of the
Stone-\u{C}ech compactification. Since $p$ is $G$-equivariant, given an
a.i.c.m.
$(m_i)_{i\in I}$ for $(X,G)$, the net of maps $y \mapsto m_{i}^{p(y)}$
defines an
a.i.c.m. for $(\beta G, G)$. \qed

\enddemo

An easy reformulation of Propositions 2.2 and 2.5 in this context, gives
the following useful
characterization.

\proclaim{4.4 Proposition}
Let $G$ be a countable discrete group. The following
conditions are equivalent:
\roster
\item $G$ is amenable at infinity.
\item There exists a sequence $(g_n)_{n\geq 1}$ of nonnegative functions on
$G\times G$ with support in a tube such that
\itemitem{\rm a)} for each $n$ and each $s$, $\sum_{t\in G} g_n(s,t) = 1$;
\itemitem{\rm b)} $\lim_n \sum_{u\in G} |g_n(s,u) - g_n(t,u)| = 0$, uniformly
on tubes.
\item There exists a sequence $(\xi_n)_{n\geq 1}$ of  functions on
$G\times G$ with support in a tube such that
\itemitem{\rm a)} for each $n$ and each $s$, $\sum_{t\in G} |\xi_n(s,t)|^2
= 1$;
\itemitem{\rm b)} $\lim_n \sum_{u\in G} |\xi_n(s,u) - \xi_n(t,u)|^2 = 0$,
uniformly
on tubes.
\item There exists a sequence $(h_n)_{n\geq 1}$ of bounded positive type
functions on
$G\times G$ with support in a tube such that
 $\lim_n h_n = 1$, uniformly on tubes.
\endroster
\endproclaim

\demo{Proof} To go from Propositions 2.2 and 2.5 to the above statement and
conversely, we just
need a simple change of variables, and we also use the universal property of
$\beta G$. Assume for instance the existence of a sequence $(h_n)$ as in
(4) above.
For $(s,t) \in G\times G$, we set $k_n(s,t) = h_n(s^{-1}, s^{-1}t)$. Now,
for every
fixed $t\in G$, we extend $s \mapsto k_n(s,t)$ to a continuous function on
$\beta G$. In this way we get a continuous function with compact support on
$\beta G\times G$. For
$x\in G$, $t_1,\dots,t_k \in G$, we have
$$ k_n(t_{i}^{-1}x, t^{-1}_{i}t_j) = h_n(x^{-1}t_i, x^{-1}t_j).$$
It follows that $(k_n)$ is a net of continuous positive type functions on
$\beta G\times G$ with compact support. Morever for every finite subset
$F\subset G$ we have
$$\aligned
\sup_{(x,t)\in \beta G\times F} |k_n(x,t) - 1| & =
\sup_{(x,t)\in G\times F} |k_n(x,t) -1|\\
&= \sup_{(x,t)\in G\times F} |h_n(x^{-1}, x^{-1}t) - 1|\\
& = \sup \{|h_n(s,t) - 1| : (s,t) \in G\times G, s^{-1}t \in F\}.
\endaligned
Ò$$
Therefore $(h_n)$ goes to one uniformly on compact subsets of $\beta
G\times G$. The
converse is even easier and the other assertions are proved in the same
way. \qed
\enddemo

\definition{4.5 Definition \cite{\bf KW}} We say that a locally compact
group $G$
is {\it exact} if for every $G$-equivariant exact sequence
$$0 \longrightarrow I \longrightarrow A \longrightarrow A/I \longrightarrow 0$$
of $G$-$C^*$-algebras, the sequence
$$0 \longrightarrow C_{r}^*(G,I)  \longrightarrow C_{r}^*(G,A)
 \longrightarrow C_{r}^*(G,A/I))  \longrightarrow 0 $$
is exact.
\enddefinition

Recall also that the corresponding sequence with the full crossed products
instead of
the reduced ones is automatically exact. For the basic facts on exact
$C^*$-algebras, we refer to \cite{\bf Wa}.

\proclaim{4.6 Theorem} Let $G$ be a locally compact group and consider the
following
conditions :
\roster
\item $G$ is amenable at infinity.
\item $G$ is exact.
\item For every exact $G$-$C^*$-algebra $B$, the crossed product
$C_{r}^*(G,B)$
is an exact $C^*$-algebra.
\item $C_{r}^*(G)$ is an exact $C^*$-algebra.
\endroster
Then (1) $\Rightarrow$ (2) $\Rightarrow$ (3) $\Rightarrow$ (4). Moreover, if
$G$ is discrete, then (4) $\Rightarrow$ (1), and therefore all these
conditions are
equivalent.

\endproclaim

\demo{Proof} (1) $\Rightarrow$ (2). Assume that $(Y,G)$ is an amenable
transformation
group with $Y$ compact. Let
$0 \rightarrow I \rightarrow A \rightarrow A/I \rightarrow 0$
be a $G$-equivariant exact sequence of $G$-$C^*$-algebras. Then
$$ 0 \rightarrow C(Y)\otimes I \rightarrow C(Y)\otimes A \rightarrow
C(Y)\otimes A/I \rightarrow 0$$
is obviously a $G$-equivariant exact sequence of $G$-$C(Y)$-algebras, when
these algebras are endowed with the diagonal $G$-actions. The corresponding
sequence
$$0 \rightarrow C^*(G,C(Y)\otimes I)  \rightarrow C^*(G,C(Y)\otimes A)
 \rightarrow C^*(G,C(Y)\otimes A/I))  \rightarrow 0 $$
of full crossed products is exact. On the other hand by Theorem 3.4, the
full and reduced
crossed products are the same. Then we obtain the following commutative
diagramm
$$
\CD
0 \longrightarrow C_{red}^*(G,I) @>{i_r} >> C_{red}^*(G,A)@>{q_r} >>
C_{red}^*(G,A/I) \longrightarrow 0\\
 @VV{\phi_I}V @VV{\phi_A}V  @VV{\phi_A}V\\
0 \rightarrow C_{red}^*(G,C(Y,I))@> >> C_{red}^*(G,C(Y,A)) @> >>
C_{red}^*(G,C(Y,A/I)) \rightarrow 0
\endCD
$$
where the bottom line is exact and the vertical arrows are the inclusions
induced
by the corresponding $G$-equivariant embeddings of $I$, $A$ and $A/I$
into $C(Y,I)$, $C(Y,A)$ and $C(Y,A/I)$ respectively.

Let us show that the first line is exact in the middle. Let $a \in
C_{red}^*(G,A)$
such that $q_r(a) = 0$. Then we have
$\phi_A(a) \in C_{red}^*(G,C(Y, I))$. Given an approximate unit $(u_\lambda)$
for $I$, the net $(\phi_I(u_{\lambda}))$ is an approximate unit of multipliers
of $C_{red}^*(G,C(Y, I))$. It follows that
$$\phi_A(a) = \lim_{\lambda} \phi_I(u_{\lambda})\phi_A(a) =\lim_{\lambda}
\phi_I(u_\lambda a) \in \phi_I(C_{r}^*(G,I)),$$
and thus $a \in C_{r}^*(G,I)$ since $\phi_A$ is injective.

(2) $\Rightarrow$ (3). Let $B$ be an exact $G$-$C^*$-algebra and
$0 \rightarrow I \rightarrow A \rightarrow A/I \rightarrow 0$
be a short exact sequence. Then
$$0 \rightarrow I\otimes B \rightarrow A \otimes B
\rightarrow A/I \otimes B \rightarrow 0$$
is an exact sequence of $G$-$C^*$-algebras where $G$ acts trivially on the
first
component of each tensor product. By (2),
$$0 \rightarrow C_{r}^*(G,I\otimes B) \rightarrow C_{r}^*(G,A \otimes B)
\rightarrow C_{r}^*(G,A/I \otimes B) \rightarrow 0$$
is an exact sequence, and we just have to observe that it coincides with
the following one:
$$0 \rightarrow I\otimes C_{r}^*(G,B) \rightarrow A \otimes C_{r}^*(G,B)
\rightarrow A/I \otimes C_{r}^*(G, B) \rightarrow 0.$$

(3) $\Rightarrow$ (4) is obvious. Assuming now that $G$ is discrete,
the proof of (4) $\Rightarrow$ (1) works as that of (4) $\Rightarrow$ (1)
in Theorem 3.4 and is in fact easier since $X$ is here reduced to a point.
By Kirchberg's characterization of exactness (see \cite{\bf Ki} or
\cite{\bf Wa}), there exists a
sequence $(\Phi_k = Q_k\circ P_k)$ of unital completely positive maps
from $C_{r}^*(G)$ into ${\Cal B}(l^2(G))$, where $P_k : C_{r}^*(G) \rightarrow
M_{n_k}(\C)$ and $Q_k :  M_{n_k}(\C) \rightarrow {\Cal B}(l^2(G))$ are
unital and completely positive, such that for every $a\in C_{r}^*(G)$,
$\lim_k \Phi_k(a) = a$.

For $(s,t) \in G\times G$, we set
$$h_k(s,t) = \langle \delta_{s^{-1}},
\Phi_k(\lambda(s^{-1}t))\delta_{t^{-1}}\rangle,$$
$(\delta_s)_{s\in G}$ being the canonical orthonormal basis of $l^2(G)$.
Obviously, $(h_k)$ is a sequence of bounded positive type kernels on
$G\times G$ which
tends to $1$ uniformly on tubes. Using Lemma 3.2 with $X$ reduced to a
point, we may
replace $(h_k)$ by a sequence of functions which, in addition, have their
support in
tubes. \qed
\enddemo

\remark{4.7 Remarks} More generally (4) implies (1) for every locally
compact group
whose unit has a basis of neighbourhoods which are invariant by inner
automorphisms.
We do not know whether (4) implies (1) for every locally compact group nor
even if
(2) implies (1).

For $G$ discrete, the equivalence between (2) and (4) is a result of
Kirchberg and
Wassermann \cite{\bf KW}, and the equivalence between conditions (4) and
(1) was also observed
by Ozawa \cite{\bf Oz}.
\endremark

Let us end these notes by a proof of the known fact \cite{\bf Yu}
that the finitely generated groups
which are amenable at infinity are uniformly embeddable into Hilbert spaces.

\definition{4.8 Definition}
A function $f: X \rightarrow Y$ between metric spaces is a {\it uniform
embedding}
if for every $R  > 0$ there exists $S > 0$ such that
$$\aligned d(x_1,x_2) \leq R & \Rightarrow d(f(x_1), f(x_2)) \leq S, \\
d(f(x_1), f(x_2)) \leq R & \Rightarrow d(x_1,x_2) \leq S.
\endaligned
$$
\enddefinition

\proclaim{4.9 Proposition} Let $G$ be a finitely generated group, endowed
with the length
function $l$ associated to the choice of a finite set of generators. Assume
that $G$ is amenable
at infinity. Then $G$, equipped with the metric defined by $l$ is uniformly
embeddable
into a Hilbert space.
\endproclaim

\demo{Proof} Consider a sequence $(\xi_n)$ of kernels on $G\times G$, with
support in
tubes such that

a) $\displaystyle \forall n\geq 1, \forall s \in G, \sum_{t\in G}
|\xi_n(s,t)|^2  = 1$,

b) $\displaystyle \forall n\geq 1, \,\sup\{ \sum_{u\in G} |\xi_n(s,u) -
\xi_n(t,u)|^2
: (s,t) \in G\times G, l(s^{-1}t) \leq n\} \leq \frac{1}{n^2}.$

For $s \in G$, $\xi_{n}^s$ will denote the element $t \mapsto \xi_n(s,t)$
of $l^2(G)$.
Let us take $H = l^2(G) \otimes l^2(\N^*)$ and consider the map
$$s\mapsto f(s) = \oplus_{n\geq 1} (\xi_{n}^s - \xi_{n}^e)$$
from $G$ into $H$. Then we have
$$\| f(s) - f(t)\|^2 = \sum_{k} \|\xi_{k}^s  - \xi_{k}^t\|_{2}^2$$
and if $l(s^{-1}t) \leq n$, we have
$$\| f(s) - f(t)\|^2 \leq 4n + \sum_{k>n} \frac{1}{k^2}.$$

Conversely, let $\psi(n)$ be an integer such that for $1\leq p \leq n$,
$$\hbox{Supp}\,\xi_p \subset \{ (s,t) : l(s^{-1}t) \leq \psi(n)\}.$$
In particular, if $l(s^{-1}t) > 2 \psi(n)$ we have
$$\hbox{Supp}\,\xi_{p}^s \cap \hbox{Supp}\,\xi_{p}^t = \emptyset
\quad \hbox{for} \quad 1\leq p\leq n,$$
and thus $\|\xi_{p}^s - \xi_{p}^t \|_{2}^2 = 2$.
It follows that if $\| f(s) - f(t)\| < \sqrt{2n}$, we have
$l(s^{-1}t) \leq 2 \psi(n)$. \qed
\enddemo

\remark{4.10 Remark} In a recent paper \cite{\bf Gro}, Gromov has announced
the existence of
finitely generated discrete groups
which are not uniformly embeddable into Hilbert spaces,
and therefore are not exact.
\endremark

\definition{4.11 Definition} Let $(X,G)$ be a locally compact transformation
group with $X$ and $G$ locally compact. We say that $(X,G)$ is {\it amenable
at infinity} if there exists an amenable transformation group $(Y,G)$ with
a proper surjective $G$-equivariant map $p : Y \rightarrow X$.
\enddefinition

In this framework, where $G$-$C^*$-algebras have to be replaced by
$G$-$C_0(X)$-algebras, we get a statement similar to that of Theorem 4.4.

\proclaim{4.12 Theorem} Let $(X,G)$ be a transformation group, and consider
the following
conditions :
\roster
\item $(X,G)$ is amenable at infinity.
\item For every $G$-$C_0(X)$-equivariant exact sequence
$$0 \longrightarrow I \longrightarrow A \longrightarrow A/I \longrightarrow 0$$
of $G$-$C_0(X)$-algebras, the sequence
$$0 \longrightarrow C_{r}^*(G,I)  \longrightarrow C_{r}^*(G,A)
 \longrightarrow C_{r}^*(G,A/I))  \longrightarrow 0 $$
is exact.
\item For every exact $G$-$C_0(X)$-algebra $B$, the crossed product
$C_{r}^*(G,B)$
is an exact $C^*$-algebra.
\item $C_{r}^*(X\croi\,G)$ is an exact $C^*$-algebra.
\endroster
Then (1) $\Rightarrow$ (2) $\Rightarrow$ (3) $\Rightarrow$ (4). Moreover, if
$G$ is discrete, then (4) $\Rightarrow$ (1), and therefore all these
conditions are
equivalent.

\endproclaim

More generally, there is a notion of locally compact groupoid
amenable at infinity. This will be the subject of a forthcoming joint work
with Jean Renault.

\heading
References
\endheading
\item{[Ad]} S. Adams: {\it Boundary amenability for word hyperbolic groups
and an
application to smooth dynamics of simple groups}, Topology, {\bf 33}
(1994), 765-783.

\item{[AD]} C.~Anantharaman-Delaroche: {\it Syst\`emes dynamiques non
commutatifs
et \linebreak moyennabilit\'e}, Math. Ann., {\bf 279} (1987), 297-315.

\item{[ADR]} C.~Anantharaman-Delaroche, J. Renault: {Amenable groupoids},
to appear, Monographie de L'Enseignement Math\'ematique, Gen\`eve.

\item{[Co]} A.~Connes: {\it Classification of injective factors},
Ann. of Math., {\bf 104} (1976), 73-115.

\item{[Di]} J.~Dixmier: Les $C^*$-alg\`ebres et leurs repr\'esentations,
Gauthier-Villars, Paris, 1964.

\item{[Fu]} H. Furstenberg: {\it Boundary theory and stochastic processes in
homogeneous spaces}, in Harmonic analysis on homogeneous spaces, Proc.
Symp. Pure and Appl., Math.
{\bf 26} (1973), 193-229.

\item{[Gre]} F.~P.~Greenleaf: Invariants means on topological groups,
Van Nostrand, New-York, 1969.

\item{[Gro]} M.~Gromov: {\it Spaces and questions}, Preprint 1999.

\item{[Hi]} N.~Higson: {\it Bivariant $K$-theory and the Novikov
conjecture}, Preprint 1999.

\item{[HR]} N.~Higson, J.~Roe: {\it Amenable group actions and the Novikov
conjecture},
Preprint 1998.

\item{[Ki]} E.~Kirchberg: {\it On non-semisplit extensions, tensor products and
exactness of group $C^*$-algebras}, Invent. Math., {\bf 112} (1993), 449-489.

\item{[KW]} E.~Kirchberg, S.~Wassermann: {\it Exact groups and continuous
bundles
of $C^*$-algebras}, Math. Ann., {\bf 315} (1999), 169-203.

\item{[La]} E.~C.~Lance: Hilbert $C^*$-modules. A toolkit for
operator algebraists, London Math. Soc. Lecture Note Series {\bf 210},
Cambridge Univ. Press, 1995.

\item{[Oz]} N. Ozawa: {\it Amenable actions and exactness for discrete groups},
Preprint 2000.

\item{[Pe]} G.~K.~Pedersen: $C^*$-algebras and their automorphism groups,
Academic Press, 1979.

\item{[Wa]} S.~Wassermann: {\it Exact $C^*$-algebras and related topics},
Lecture Notes Series, {\bf 19}, Global Analysis Research Center, Seoul
National University,
1994.

\item{[Yu]} G.~Yu: {\it The Baumes-Connes conjecture for spaces which admit
a Uniform embedding into Hilbert space}, Invent. Math., to appear.

\enddocument